\documentclass[preprint,12pt]{elsarticle}

\usepackage{amssymb}
\usepackage{amsthm}

\newcommand{\sysn}{\left\{\begin{array}{rcl}}
\newcommand{\sysk}{\end{array}\right.}

\newtheorem{theorem}{Theorem}[section]

\theoremstyle{example}

\theoremstyle{definition}

\begin{document}

\begin{frontmatter}


%%\title{Topological-Algebraic Properties of Function Space with Set-Open Topology\tnoteref{label1}}
%%\tnotetext[label1]{}
%%\author{Alexander V. Osipov\corref{cor1}\fnref{label2}}
%%\ead{OAB@list.ru}

\title{Remarks on the Ostrovsky's Theorem  \tnoteref{label1}}


\author{Alexander V. Osipov}

\ead{OAB@list.ru}

%\tnotetext[label1]{The research has been supported by .}

\address{Krasovskii Institute of Mathematics and Mechanics, Ural Federal
 University,

 Ural State University of Economics, Yekaterinburg, Russia}

\begin{abstract}

In this paper we prove that the condition {\it 'one-to-one'} of
continuous open-resolvable mapping is necessary in the Ostrovsky's
Theorem (Theorem 1 in \cite{Os1}). Also we get that the
Ostrovsky's Problem (\cite{Os2}, Problem 2) ({\it Is every
open-$LC_n$ function between Polish spaces piecewise open for
$n=2,3,...$ ?}) has a negative solution for any $n>1$.


\end{abstract}


\end{frontmatter}

% \linenumbers


\section{Introduction}

In the following definitions we will suppose that $X$ is a
subspace of the Cantor set $\mathbf{C}$.

A function $f: X\mapsto Y$ is called {\it piecewise open} if $X$
admits a countable, closed and disjoint cover $\mathcal{V}$, such
that for each $V\in \mathcal{V}$ the restriction $f|V$ is open.

Recall, that a subset $E$ of a metric space $X$ is {\it
resolvable} \cite{kr}, if for each nonempty closed in $X$ subset
$F$ we have $cl_{X}(F\cap E)\cap cl_{X}(F\setminus E)\neq F$.

If $E\subset X$ is resolvable, then $E$ is $\Delta^0_2$-set in $X$
and vice versa if the space $X$ is Polish.

Recall that a subsets $X$ is $LC_n$-set if it can be written as
union of $n$ locally closed in $X$ sets (a set is locally closed
if it is the intersection of an open set and a closed set.) Every
$LC_n$-set (constructible) set is resolvable.

A mapping $f$ is open if it maps open sets onto open ones. More
generally, for $n\in \omega$ a mapping $f$ is said to be {\it
open-resolvable} (open-$LC_n$) if $f$ maps open set onto
resolvable ($LC_n$-set) ones.

In the following definitions we will suppose that $X$ is a
subspace of the Cantor set $\mathbf{C}$.

A piecewise open function $f:X\mapsto Y$ is called {\it
scatteredly open} if, in addition, the cover $\mathcal{V}$ is
scattered, that is: for every nonempty subfamily
$\mathcal{T}\subset \mathcal{V}$ there is a clopen set $G\subset
X$ such that $\mathcal{T}_G=\{T\in \mathcal{T}: T\subset G\}$ is a
singleton and $T\bigcap G=\emptyset$ for every $T\in
\mathcal{T}\setminus \mathcal{T}_G$.

\section{Main result}

A.V. Ostrovsky proved the interesting results

\begin{theorem}\label{th1}(Theorem 1 in \cite{Os1}) Let $X$ be subspace of
the Cantor set $\mathbf{C}$, and $f: X \mapsto Y$ a continuous
bijection. If the image under $f$ of every open set in $X$ is
resolvable in $Y$, then $f$ is scatteredly open, and, hence, $f$
is scattered homeomorphism.

\end{theorem}

\begin{theorem}\label{th2}(Proposition 3.2 in \cite{Os3}) Every
continuous open-$LC_1$ function $X\mapsto Y$ onto a metrizable
crowded space $Y$ is open.
\end{theorem}

\medskip

In (\cite{Os2}, Problem 2) A.V. Ostrovsky posed the following
\medskip

{\bf Problem.} Is every open-$LC_n$ function between Polish spaces
piecewise open for $n=2,3,... ? $

\medskip

We prove that

$\bullet$ the condition {\it 'one-to-one'} of mapping $f$ in
Theorem \ref{th1} is necessary.

$\bullet$ the Ostrovsky's Problem has a negative solution for
$n=2$ (hence for every $n>1$).

\medskip

{\bf Example.} Let $\mathbf{C}$ be the Cantor set such that
$\mathbf{C}\subset [0,1]$. As usually, we starts be deleting the
open middle third $(\frac{1}{3}, \frac{2}{3})$ from the interval
$[0,1]$, leaving two segments: $P_1=C_0\cup
C_2=[0,\frac{1}{3}]\cup [\frac{2}{3},1]$. Next, the open middle
third of each of these remaining segments is deleted, leaving four
segments: $P_2=C_{00}\cup C_{02}\cup C_{20}\cup
C_{22}=[0,\frac{1}{9}]\cup
[\frac{2}{9},\frac{1}{3}]\cup[\frac{2}{3},\frac{7}{9}]\cup
[\frac{8}{9},1]$. This process is continued ad infinitum, where
the $n$th set is $P_n=\frac{P_{n-1}}{3}\cup (\frac{2}{3}+
\frac{P_{n-1}}{3})$ for $n\geq 1$, and $P_0=[0,1]$.

The Cantor ternary set contains all point in the interval $[0,1]$
that are not deleted at any step in this infinite process:

$\mathbf{C}:=\bigcap\limits_{n=1}^{\infty} P_n$.

 Let us fix a countable dense set
$\{(a_n, b_n) : n\in \omega \}$ in $\mathbf{C}\times \mathbf{C}$
such that $a_n\neq a_m$ and $b_n\neq b_m$, for $n\neq m$, for each
$n$ pick $a_{n,i}\mapsto a_n$ such that $a_{n,i}\neq a_m$,
$a_{n,i}\neq a_{m,j}$ for $(n,i)\neq (m,j)$, and $|a_{n,i}-a_n|<
\frac{1}{n}$.

Consider the standard clopen base $\mathcal{B}:=\{\mathbf{C}\cap
C_{s_1,...,s_k}: s_i\in \{0,2\}, i\in\overline{1,k}$,
$k\in\omega\}$ in $\mathbf{C}$,  and we enumerate
$\mathcal{B}=\{B_n : n\in \omega\}$ such that $b_n\in B_n$ for
every $n\in \omega$.

Let $X=(\mathbf{C}\times \mathbf{C})\setminus \bigcup\limits_{n,i}
\{a_{n,i}\}\times B_n$. Note that $X$ is $G_\delta$-set of
$\mathbf{C}\times \mathbf{C}$. It follows that $X$ is
$\check{C}$ech-complete and, moreover, it is Polish space.

Let $\pi | X: X \mapsto \mathbf{C}$ be the restriction to $X$ of
the projection $\pi : \mathbf{C}\times \mathbf{C} \mapsto
\mathbf{C}$ onto the first coordinate. Note that
$\pi(X)=\mathbf{C}$ because $diam \mathbf{C}> diam B_n$ for any
$n\in \omega$.

Suppose $X=\bigcup\limits_{n\in \omega} X_n$ is a countable union
of closed subsets $X_n$ (apply the Baire Category Theorem), there
is $X_{m}$ such that $V=Int X_{m}\neq \emptyset$.

Since the set $\{(a_n, b_n) : n\in \omega \}$ is dense in $X$,
there are  $n'\in \omega$ and $W\in \mathcal{B}$ such that a point
$(a_{n'},b_{n'})\in ((W\times B_{n'})\cap X)\subset V$. Since the
set $\{(a_n, b_n) : n\in \omega \}$ is dense in $(W\times
B_{n'})\cap X$, choose $n''\in \omega$ such that $n''>n'$, $2*diam
B_{n''}<diam B_{n'}$ and $(a_{n''},b_{n''})\in ((W\times
B_{n''})\cap X)\subset (W\times B_{n'})\cap X$.  Then $\pi | X_{m}
: X_{m} \mapsto \pi(X_{m})$ is not open at $(a_{n''}, b_{n''})$
because of $\pi((W\times B_{n''})\cap X)$  is not contains
$\{a_{n'',i}: i\in \omega\}$ and, hence, it is not open set of
$\pi(X_{m})$. Therefore $\pi | X$ is not piecewise open and,
hence, is not scatteredly open.

Let $U\subset \mathbf{C}\times \mathbf{C}$ be open. We have to
check that $\pi(U\cap X)\in \Delta^0_2$.

Construct for every point $(a,b)\in U\cap X$ a sets $W(a)$ and
$B(b)$ such that

$\bullet$  $a\in W(a)\in \mathcal{B}$, $b\in B(b)\in \mathcal{B}$
and $(W(a)\times B(b))\cap X\subset U$.

$\bullet$  if $a\neq a_m$ for any $m\in \omega$, then
$\pi((W(a)\times B(b))\bigcap X)=W(a)$.

$\bullet$ if $a=a_m$ for some $m\in \omega$, then $\pi((W(a)\times
B(b))\bigcap X)=W(a)\setminus \{ a_{m,i_j}: j\in \omega \}$ for
some subsequence $\{ a_{m,i_j}: j\in \omega \}\subseteq \{
a_{m,i}: i\in \omega \}$.

 Case 1. Let $a\neq a_m$ for any $m\in \omega$, one can
choose $W$, $B(b)=B_{n'}\in \mathcal{B}$ such that $a\in W$, $b\in
B(b)$, $(W\times B(b))\cap X\subset U$ and $B(b)\setminus B_n\neq
\emptyset$ for all $n>n'$. Since $a\neq a_m$ for any $m\in
\omega$, then there exist $W(a)\in \mathcal{B}$ such that $a\in
W(a)\subset W$ and $W(a)\cap \{a_i\cup \{a_{i,j}: j\in \omega\} :
i\in \overline{1,n'}\}=\emptyset$. Then $\pi((W(a)\times
B(b))\bigcap X)=W(a)$.

Case 2. Let $a=a_m$ for some $m\in \omega$, analogically to Case
1, we can choose $B(b)\in \mathcal{B}$ such that $B(b)\setminus
B_n\neq \emptyset$ for all $n>n'>m$, and $W(a)\in \mathcal{B}$ can
choose such that $W(a)\cap \{a_i\cup \{a_{i,j}: j\in \omega\} :
i\in \overline{1,n'}$ and $i\neq m\}=\emptyset$.

 Then $W(a)\setminus \pi((W(a)\times B(b))\bigcap X)\subset \{
a_{m,i}: i\in \omega \}$, hence $\pi((W(a)\times B(b))\bigcap
X)=W(a)\setminus \{ a_{m,i_j}: j\in \omega \}=W_a\cup \{a_m\}$
where $W_a=W(a)\setminus (\{a_m\}\cup \{ a_{m,i_j}: j\in \omega
\})$ is an open in $\mathbf{C}$.

Thus $\pi(U\cap X)=\bigcup\limits_{(a,b)\in U\cap X}
\pi((W(a)\times B(b))\bigcap X)=$

$=(\bigcup\limits_{(a,b)\in U\cap X, a\neq a_m} W(a))\cup
(\bigcup\limits_{(a,b)\in U\cap X, a=a_m} W_a \cup \{a_m\})$.

 By definition of the clopen base $\mathcal{B}$,
$\pi(U\cap X)=S\cup D$ where $S=(\bigcup\limits_{(a,b)\in U\cap X,
a\neq a_m} W(a))\cup (\bigcup\limits_{(a,b)\in U\cap X, a=a_m}
W_a)$ is an open set in $\mathbf{C}$ and $D=\{a_{m_k}: k\in \omega
\}$ is a discrete in itself such that $S\bigcap D=\emptyset$.
Indeed, by Case 2, for every $a_{m_k}\in D$ there is
$W(a_{m_k})\in \mathcal{B}$ such that $a_{m_k}\in W(a_{m_k})$ and
$W(a_{m_k})\bigcap \{a_{m_i} : i\in \omega, i\neq k \}=\emptyset$.
It follows that $D$ is a discrete in itself, and, hence,
$\pi(U\cap X)$ is $\Delta^0_2$. Since $\pi(X)=\mathbf{C}$ is
Polish, the mapping $\pi|X$ is continuous open-resolvable.

Note that $\pi(U\cap X)=S\cup ((\bigcup\limits_{(a,b)\in U\cap X}
W(a))\bigcap \overline{D})$. It follows that $\pi(U\cap X)$ is
$LC_2$-set and, hence, $\pi|X$ is open-$LC_2$.

\smallskip

% {\bf Acknowledgment.}



%% References


\bibliographystyle{model1a-num-names}
\bibliography{<your-bib-database>}


%%\bibliographystyle{plain}

\end{document}